\documentclass[11pt]{article}

\usepackage{amsmath}
\usepackage{amssymb}
\usepackage{latexsym}
\usepackage{hyperref}
\usepackage{graphicx}
\usepackage{amsmath,amssymb,float}
\usepackage{epsfig}
\usepackage{amsthm}
\usepackage{comment}
\usepackage{pgfplots}
 \usepackage{caption}
\usepackage{tikz}
\usepackage{subcaption}     
\usepackage{enumitem}
\usepackage{algorithm}
\usepackage{algpseudocode}
\usepackage{authblk}

\newcommand{\argmax}{\operatornamewithlimits{argmax}}
\def\be{\begin{eqnarray}}
\def\ee{\end{eqnarray}}

\def\b*{\begin{eqnarray*}}
\def\e*{\end{eqnarray*}}

\newcommand{\Cn}[1]{\text{\rm C}^{#1}}

\newtheorem{Theorem}{Theorem}[part]
\newtheorem{Definition}[Theorem]{Definition}
\newtheorem{Proposition}[Theorem]{Proposition}

\newtheorem{Lemma}[Theorem]{Lemma}
\newtheorem{Corollary}[Theorem]{Corollary}
\newtheorem{Remark}[Theorem]{Remark}

\makeatletter \@addtoreset{equation}{section}

\@addtoreset{Theorem}{section}

\def\diff{\text{d}}

\def \E{\mathbb{E}}
\def \F{\mathbb{F}}

\def \M{\mathbb{M}}

\def \P{\mathbb{P}}
\def \Q{\mathbb{Q}}
\def \R{\mathbb{R}}

\def\Ac{{\mathcal A}}
\def\Bc{{\mathcal B}}

\def\Fc{{\mathcal F}}

\def\Lc{{\mathcal L}}

\def\Qc{{\mathcal Q}}

\def \Int{\displaystyle\int}

\newtheorem{innercustomthm}{Assumption}
\newenvironment{assp}[1]
  {\renewcommand\theinnercustomthm{#1}\innercustomthm}
  {\endinnercustomthm}

\def\={\;=\;}
\def\.{\;.}

\def\eps{\varepsilon}

 \def\normeL2#1{\left\|{#1}\right\|_{L^2}}
\newcommand {\Chi} {{\bf \raise 1.5pt \hbox{$\chi$}}}
 \title{Impact of carbon market on production emissions}

\author[1]{Arash {\sc Fahim}}
\affil[1]{Florida State University}
\author[2]{Nizar {\sc Touzi}}
\affil[2]{Ecole Polytechnique}

\date{September 5, 2015}

\begin{document}

\maketitle

\begin{abstract}
The aim of this paper is to address the effect of the carbon emission allowance market on the production policy of a large polluter production firm. We investigate this effect in two cases; when the  large polluter cannot affect the risk premium of the allowance market, and when it can change the risk premium by its production. In this simple model, we ignore any possible investment of the firm in pollution reducing technologies. We formulate the problem of optimal production by a stochastic optimization problem. Then, we show that, as expected, the market reduces the optimal production policy in the first case if the firm is not given a generous initial cheap allowance package. However, when the large producer activities can change the market risk premium, the cut on the production and consequently pollution cannot be guaranteed. In fact, there are cases in this model when the optimal production is {\it always} larger than expected, and an increase in production, and thus pollution, can increase the profit of the firm. We conclude that some of the  parameters of the market which contribute to this effect can be wisely controlled by the regulators in order to diminish this manipulative behavior of the firm.


\end{abstract}

%%%%%%%%%%%%%%%%%%%%%%%%%%%%%%%%%%%%%%%%%%%%%%%%%%%%
%%%%%%%%%%%%%%%%%%%%%%%%%%%%%%%%%%%%%%%%%%%%%%%%%%%%
%%%%%%%%%%%%%%%%%%%%%%%%%%%%%%%%%%%%%%%%%%%%%%%%%%%%
\section{Introduction}
\label{Sec:int}
The long term costs of global warming is believed to be significantly more than the cost of controlling it by reducing the pollution due to greenhouse gases (see \cite{mccarthy01}). The Kyoto protocol in 1997 concerns with the reduction of the greenhouse gases including $\text{CO}_2$ and is accepted by several nations. These nations agreed to set goal on reduction of greenhouse gases and implement plans to reach  the  goals. One of the popular plans in the so-called {\it cap-and-trade} scheme adopted by several nations including members of European Union. The principle behind standard cap-and-trade is simple. Regulator marks the polluter installations and set a cap on the total emission at the end of a specific period of time and issue a number of  allowances  equal  to the cap. Then, they allocate the allowances to the those installations.  If the cap is reached, then all installations are mandated to pay a predetermined penalty or present sufficient allowances. At the same time, there is a market where they can trade for the allowances; if a firm does not need its allowances for whatever reason, they can sell it to those who want to produce more. Ther allowance papers are worthless, if the cap is not reached. 

There are several variations of standard cap-and-trade scheme, running around the world, e.g. European Union Emission Trading Scheme (EU ETS), US REgenial CLean Air Incentive Market (RECLAIM) or Regional Greenhouse Gas Incentive (RGGI). Although based on the same principle, they may differ in certain details. Some of the specifics of the cap-and-trade market is to make it work more efficiently. For instance to avoid a sharp drop in the allowance price when the cap is not reachable, they can store their allowances for the next period of the scheme by paying an extra fee per allowance, which is referred to as {\it banking}. Also, regulator can distribute initial allowances to involved installations  for free, or they can set an auction to distribute them, or a mixture of both. For example in the third phase of EU ETS  2013-2020, 40\% of the allowances are distributed by auctioning, while in the first two phases 2006-2008 and 2008-2012, the allowances are distributed for free. Also, the cap can be set regionally where each region has its own cap on the total emission of the region, or globally where all regions have one cap on the total emission. The later creates a less stressful market. The efficiency of the design of the market is  comprehensively  studied in \cite{cfhp10} where they show in the standard cap-and-trade system the presence of windfall profit based on the real market data and propose a more efficient allocation scheme. In \cite{cfh-2-09}, the authors study some alternatives to the standard cap-and-trade system which can potentially lead to less windfall profit for a dominant player in the market and less cost for the consumers of the product of the polluter firms. More precisely, in addition to initial auctioning, they propose to allocate  part of the allowance over time, which make it more flexible for the regulator to achieve its pollution reduction target with less cost on the economy.

Several studies target the dynamics of the price of the allowance.  In the presence of banking, the price of allowance in the current period can be viewed as an option on the price of allowance for the next period (See for example \cite{cv09}). This approach is obviously not capable of explaining the price of the most farthest period, and therefore, it is important to take a different approach  to explain the dynamics for the price in the last period. In \cite{cfh-1-09}, by adopting a stochastic game setting in discrete-time, the authors show the existence of a Nash equilibrium for the price of emission allowance, which appears to be a martingale. In addition, they show that this equilibrium price is equal to the marginal cost of total abatement which can be obtained by solving the problem of minimizing the total abatement cost in the market. In \cite{cdet13}, they study the formation of equilibrium price in a continuous-time setting through a system of forward-backward SDEs (FBSDE) with singularity at terminal time. In their study, the singularity of the terminal condition is two-folded; one caused by the discontinuity of terminal condition of backward equation and the other by the degeneracy of the forward equation at terminal time. Then for some specific pathological examples, they showed the existence and the uniqueness of the solution to the system of FBSDEs, by approximating the terminal condition with smooth functions in a certain manner. Beside  modeling allowance price by FBSDE, their main contribution is to bold the difficulties caused by discontinuous terminal condition of backward component and the demand for a more inclusive theory of FBSDEs in this case.

In this paper, we analyze the effect of the cap-and-trade scheme in reducing the carbon emission through reduction on the optimal production of the relevant production firms. The setting of this paper is similar to the third phase of EU ETS where one EU-wide cap on total emission is imposed. The firm's objective is to maximize its utility on wealth, which is made of the profit gained from production and the value of  carbon allowance portfolio, over its production policy and its portfolio strategy. Via standard duality, we manage to first solve  the utility maximization problem over the portfolio strategy for a fixed production strategy. Then, we manage to derive a stochastic control problem over the optimal portfolio strategy only. In the Markovian case, the stochastic control problem for optimal production can be handled by the Hamilton-Jacobi-Bellman (HJB) equation. The terminal condition of the HJB equation is discontinuous, which causes a  similar challenge as in \cite{cdet13}. However, we partially avoid this challenge by assuming that the forward equation modeling the state process is non-degenerate and therefore has no atom at the discontinuity of the terminal condition. 

We further categorize the relevant firms by their impact on the risk premium of the price of carbon allowance. We first present the model for a  small producer  which is a price taker and can neither change the allowance price nor the total  emission significantly. We use the change in the production of  a small firm as a benchmark to make the comparison. We observe that the market always reduces the optimal production policy of a large firm who cannot affect the risk premium of the allowance price, if it is not given too much of free initial allowances. However, this study shows that a large producer with impact on the risk premium can take a manipulating role and its optimal production behavior does not necessarily reduce the emission. The key result to establish this comparison is that the price of the carbon allowance is equal to negative sensitivity of the value function of the firm with respect to its emission.

To better address the manipulative nature of a large producer, we make certain simplification in this model.  We assume that the  allowance market is complete. In addition, we ignore the effect of abatement and leave it for the future research. As for the profit function of the firm, it usually depends on the price of raw material and the product of the producer which are all stochastic. Here, we eventually ignore the stochastic nature of the profit of the firm in establishing the main results. This can be justified for a period of time in which the supply and demand remain stationary and non-volatile.
We also assume that the emission dynamics is governed by an It\^o process. This process usually represents the perception of the firm on the total emission. The total emission is not revealed by the regulators until the end of the period. While in the standard cap-and-trade each inclusive firm has to adjust its position on allowances yearly, here we assumed that the firm only needs to provide sufficient allowance only at the end of the period. % for the sake of simplicity, we do not extend this study to multi period with banking. 

%Although the idea behind all these markets is cap-and-trade, they are different in certain details. A standard cap-and-trade system allocates allowances at the beginning of the period for no charge based on regulator's discretion and penalizes the excess emission units of each installation at the end of the period or even shorter sub-periods. Alternatively, the regulator can distribute the allowances by auctioning at the beginning or even over the duration of the period. Banking also offers to transfer the allowance over the future periods by paying a certain charge. Also, it can impose the penalties only if the total emission exceeds a global cap, as it is so in the third phase of EU ETS.  About the efficiency of the design of the market, \cite{cfhp10} comprehensively studies the standard cap-and-trade system and verifies the presence of windfall profit based on the real market data and proposes a more efficient allocation scheme. In \cite{cfh-2-09}, the authors study some alternatives to the standard cap-and-trade system which can potentially lead to less windfall profit for a dominant player in the market and less cost for the consumers of the product of the polluter firms. More precisely, in addition to auctioning, they propose that the distribution of allowance over time rather than the beginning of the period will make it more flexible for the regulator to achieve its pollution reduction target with less cost on the economy. 

The study of a large producer manipulating the emission allowance market is back to \cite{me89} where they consider a monopoly (or monopsony) firm whose production impacts both allowance price and product price. They show, in a static context, that the monopoly firm can manipulate the market by  transferring the abatement costs to its rivals and as a consequence increasing the cost of production for fringe firms. In their study, the monopoly firm strives to maximize its profit subject to demand and price impact constraints. In \cite{lm06} and \cite{lm08}, the authors study the market power of a large producer in analogy to the market power in the context of exhaustible resource market in a dynamic setting where the firm decides how to buy/sell the allowance permits and how to use them over time. They show that  the large producer firm covers its total emission in a competitive manner unless the initial allocation is not sufficient for its optimal production plan. In all the above mentioned literature, the stochastic nature of the allowance price in the market is ignored.

This paper is organized as follows. In Section \ref{Sec:small}, we present a general  model and derive a characterization for the optimal production policy of a small production firm. The tool we use in this section is convex duality for utility maximization which helps us separate the trading activity of the producer and its profit from the production. In Section \ref{sec:large}, we repeat the analysis of Section \ref{Sec:small} for a large producer in two cases based on the power of impact of the large producer on the risk premium of the market. We start Section \ref{sec:large} in a general framework by using convex duality in a similar fashion to Section \ref{Sec:small}. Then later in this section, we narrow this study to Markovian setting and derive  HJB equation for the profit function of the firm. We use this HJB equation to study the impact of the large producer both in the analytical and numerical results.
In Section \ref{Sec:numerics}, we present out numerical results. Appendix sections  cover the existence and uniqueness for the HJB equation, the existence of optimal production policy, and the implementation details of the numerical results.
%%%%%%%%%%%%%%%%%%%%%%%%%%%%%%%%%%%%%%%%%%%%%%%%%%%%%%%%%%%%%%%%%
%%%%%%%%%%%%%%%%%%%%%%%%%%%%%%%%%%%%%%%%%%%%%%%%%%%%%%%%%%%%%%%%%
%%%%%%%%%%%%%%%%%%%%%%%%%%%%%%%%%%%%%%%%%%%%%%%%%%%%%%%%%%%%%%%%%
\section{Small producer with one-period carbon emission market}
\label{Sec:small}
In this section, we consider $(\Omega,\Fc,\F=\{\Fc_t,t\ge 0\},\P)$ be a filtered probability space satisfying the usual conditions which hosts one-dimensional Brownian motion $W$, and we denote by $\E_t:=\E[\cdot|\Fc_t]$ the conditional expectation operator given $\Fc_t$. 
We also consider a production firm with risk preference described by the 
utility function $U:\R\longrightarrow\R\cup\{\infty\}$ assumed to be strictly increasing, strictly concave and $\Cn{1}$ over $\{U<\infty\}$. We denote by $\pi_t(\omega,q)$ the (random) rate of profit of the 
firm for a production rate $q$ at time $t$. Here $\pi:[0,T]\times\Omega\times\R_+\to\R$ is an $\F-$progressively measurable map\footnote{An $\F-$progressively measurable map is usually defined for a mapping $\pi$ from $[0,T]\times\Omega$ to $\R$. However, we can simply extend it by calling a  mappings $\pi:[0,T]\times\Omega\times \R_+\to\R$ $\F$-progressively measurable if and only if $\pi:[0,t]\times\Omega\times \R_+\to\R$ is $\Bc([0,t])\otimes\Fc_t\otimes\Bc(\R_+)$-measurable for all $t\in[0,T]$. In this manner, if $q_\cdot(\omega)$ is a $\F$-progressively measurable process in the usual sense and $\pi$ is $\F$- progressively measurable in the extended sense, then $\pi(\cdot,\omega,q_\cdot(\omega))$ is  $\F$-progressively measurable in the usual sense.}. We shall omit $\omega$ from the notations wherever appropriate. For fixed $(t,\omega)$, we 
assume that the function $\pi_t(\cdot):=\pi(t,\cdot)$ is  $\Cn{1}$ and strictly concave in $q$, and satisfies
 \b*
0< \pi_t'(0+) &\mbox{and}& \pi_t'(\infty)<0.
 \e*
Let us denote by $\eta_t(q)$ the rate of carbon emission generated by a production rate $q$. Here, $\eta:[0,T]\times\Omega\times\R_+$ is an $\F-$progressively measurable map such that for each $(t,\omega)$, $\eta_t(q):=\eta(t,\cdot)$ is  $\Cn{1}$ and increasing in $q\in\R_+$. Then the total amount of carbon emission induced by a production policy $\{q_t,t\in[0,T]\}$ is given by 
 \be\label{defEq}
 E^q_T:=\Int_0^T \eta_t(q_t)dt.
 \ee
The aim of the carbon emission market is to incur some cost to the producer so as to obtain an overall reduction on the carbon emission.

From now on, we analyze the effect of the presence of the carbon emission market within the {cap-and-trade} scheme.

In order to model the allowance price, we introduce a state variable $Y$ given by the dynamics:
 \be\label{defY}
 dY_t
 &=&
 \mu_t dt + \gamma_t dW_t,
 \ee
where $\mu$ and $\gamma$ are two bounded $\F-$adapted processes and $\gamma>0$. 
\begin{Remark}\label{rem:Y stochastic}
The state variable $Y$ should be interpreted as the perception of the firm on the total carbon emission. Since the total emission  is only revealed at the end of the period, the process $Y$ involves uncertainty and is considered stochastic.% because of the following reasons.
%First, it depends on the total production in the environment which makes it impossible for a single firm to precisely apperceive it. Second, the amount of production and thus the emission  depends on several factors such as fuel price, weather changes, economic growth, etc; each of which are believed to be stochastic. Later on, we will assume that  a large  producer firm can benefit from knowing its contribution in the trend of   state variable.
\end{Remark}

We assume that there is one single period $[0,T]$ during which the carbon emission market is in place.  At each time $t\ge 0$, the random variable $Y_t$ indicates the  aggregated market opinion on the cumulated carbon emission. At time $T$, $Y_T\ge \kappa$ (resp. $Y_T<\kappa$) means that the cumulated total emission have (resp. not) exceeded the cap $\kappa$, fixed by the trading scheme. We simply take $\kappa=0$. Let $\alpha$ be the penalty per unit (tonne) of carbon emission. Then, the value of the carbon emission contract at time $T$ is:
 \b*
 S_T &:=& \alpha {1}_{\{Y_T\ge0\}}.
 \e*
The carbon emission allowance can be viewed as a derivative security on $Y$ defined by the above payoff. (See \cite{smw08} and \cite{ct12}.) The carbon emission market allows for trading this contract in continuous-time throughout the time period $[0,T]$. Assuming that the market is frictionless, it follows from the classical no-arbitrage valuation theory that the price of the carbon emission contract at each time $t$ is given by
 \be\label{S-1period}
 S_t
 \;:=\; 
 \E_t^{\Q}\left[S_T\right]
 \;=\;
 \alpha \Q_t\left[Y_T\ge 0\right],
 \ee
where $\Q$ is a probability measure equivalent to $\P$, the so-called equivalent martingale measure, $\E_t^{\Q}$ and $\Q_t$ denote the conditional expectation and probability given $\Fc_t$. Given market prices of the carbon allowances, the risk-neutral measure may be inferred from the market prices. 

In the present context, and in contrast to a standard taxation (Remark \ref{taxation}), production firms have more incentive to reduce emission as they have the possibility to sell their allowances on the emission market.

We now formulate the objective function of the firm in the presence of the emission market. The primary activity of the firm is the production modeled by the rate $q_t$ at time $t$. This generates a gain $\pi_t(q_t)$. The resulting carbon emission are given by $\eta_t(q_t)$. Given that the price of the allowance is available on the market, the profit on the time interval $[0,T]$ is given by $\int_0^T \pi_t(q_t)dt - S_T\Bigl(\int_0^T \eta_t(q_t)dt-E^{{(0)}}\Bigr)$,
 where $E^{(0)}$ is the number of free allowances of the firm. 
In addition to the production activity, the firm trades continuously on the carbon emission market. Let $\{\theta_t,t\ge 0\}$ be an $\F$-adapted process such that $\int_0^T \theta^2_t d\langle S\rangle_t<\infty$ $\P$-a.s.. For every $t\ge 0$, $\theta_t$ indicates the number of  allowances held by the firm at time $t$ incuding those given to the firm at time $0$, i.e. $E^{(0)}$. Under the self-financing condition, the wealth accumulated by trading in the emission market is $ x + \int_0^T \theta_tdS_t$,
where $x$ is the initial capital of the firm, including the market value of its initial allowances, i.e. $S_0E^{(0)}$.  Therefore, the total wealth of the firm at time $T$ is given by $W_T^{\theta,q}:=X_T^{x,\theta+E^{(0)}}+B^q_T$
where for all $ t\in[0,T]$
\[
X_t^{x,\theta}:=x+\int_0^t \theta_s dS_s,\;B_t^q:=\int_0^t \pi_t(q_t)dt-S_t E^q_t\;\text{ and }\;E_t^q :=\int_0^t \eta_s(q_s)ds
\]
We assume that the firm is allowed to trade with no constraint. Then, the objective function of the manager is:
 \be\label{prob:small}
 V^{(1)}
 &:=&
 \sup\left\{\E\left[U\left(W_T^{\theta,q}\right)\right]:~\theta\in\Ac,~q\in\Qc \right\},
 \ee
where $\Ac$ is the collection of all $\F-$progressively measurable processes $\{\theta_t\}_{t\ge0}$ such that $\int_0^T \theta^2_t d\langle S\rangle_t<\infty$ $\P$-a.s. and $X^{x,\theta}_t:=x+\int_0^t \theta_s dS_s$ is  bounded from below by a martingale, and $\Qc$ is the collection of all non-negative $\F-$progressively measurable bounded processes $\{q_t\}_{t\ge0}$ such that $E_T^q<\infty$. 
%{\Large \color{red} Restricted admissible strategies for Girsanov theorem. }

Notice that the stochastic integrals with respect to $S$ can be collected together in the expression of $W_T^{\theta,q}$. Since $\Ac$ is a linear subspace, it follows that the maximization with respect to $q$ and $\theta$ are completely decoupled, this problem is easily solved by optimizing successively with respect to $q$ and $\theta$. 
\begin{Proposition}\label{prop:small_theta_q}
Under the assumptions enforced in this section, the optimal production policy is independent of the utility function of the producer $U$, and obtained by solving 
\begin{equation}\label{prob:max_q1}
\sup_{q_\cdot\in\Qc} \E^\Q\left[ B^q_T \right]
\end{equation}
 where $\Q$ is the martingale measure.
Moreover,  optimal production policy $q^{(1)}$ and  optimal investment strategy $\theta^{(1)}$ are characterized by
 \begin{equation}\label{eqn:optimal_theta}
 X_T^{x,\theta^{(1)}}+B_T^{q^{(1)}}=(U')^{-1}\left(y^{q^{(1)}}\frac{d\Q}{d\P}\right),
 ~~
 x+\E^\Q\left[B_T^{q^{(1)}}\right]
 =
 \E^\Q\left[(U')^{-1}\left(y^{q^{(1)}}\frac{d\Q}{d\P}\right)\right].
 \end{equation}
\end{Proposition}
\begin{proof}  We first fix some production strategy $q$. Since the market is complete, the partial maximization with respect to $\theta$ can be performed by the duality method in \cite[Theorem 3.1]{CSW01} to obtain
 \begin{equation}\label{eqn:optimal_theta1}
 X_T^{x,\theta^{q}}+B_T^{q}=(U')^{-1}\left(y^{q}\frac{d\Q}{d\P}\right),
   x+\E^\Q\left[B^q_T\right]=
\E^\Q\left[(U')^{-1}\left(y^q\frac{d\Q}{d\P}\right)\right].
\end{equation}
Thus, problem \eqref{prob:small} can be written as
\[
 \sup_{q_.\ge 0} \E \left[U\circ(U')^{-1}\left(y^q\frac{d\Q}{d\P}\right)\right].
 \]
Notice that $U\circ(U')^{-1}$ is decreasing and the density $\frac{d\Q}{d\P}>0$. Then, finding the maximizer of the above problem  can equivalently found by solving $\inf\left\{y^q:~q_\cdot\ge 0\right\}$.
Since $(U')^{-1}$ is also decreasing, one can use \eqref{eqn:optimal_theta1} again  pass to the equivalent problem
$\sup\left\{ \E^\Q\left[ B^q_T \right]:~q_\cdot\in\Qc\right\}$
which characterizes the optimal strategy $q^{(1)}$. Finally,  the optimal investment policy is characterized by \eqref{eqn:optimal_theta} for $q=q^{(1)}$.
\end{proof}      
By using integration by parts, we can write 
\[
B_t^q=\int_0^t \bigl(\pi_s(q_s)dt-S_s \eta_s(q_s)\bigr)ds-\int_0^tE_s^qdS_s.
\]
Since  $\E^\Q[\int_0^tE_s^qdS_s]=0$, we obtain $\E^\Q[ B^q_T ]=\E^Q[\int_0^t \bigl(\pi_s(q_s)dt-S_s \eta_s(q_s)\bigr)ds]$
Problem \eqref{prob:max_q1} provides an optimal production level $q^{(1)}$ defined by the first order condition:
 \be\label{defn:q^1}
 \frac{\partial \pi_t}{\partial q}(q^{(1)}_t)
 &=&
 S_t\frac{\partial \eta_t}{\partial q}(q^{(1)}_t).
 \ee

Because of the assumptions on $\pi_t(\cdot)$  and $\eta_t(\cdot)$, we deduce immediately that $q^{(1)}_t$ is less than the {\it business-as-usual} optimal production $q^{\rm bau}$ of the firm in the absence of any restriction on the emission, which is determined by $(\partial \pi_t/\partial q)(q^{\rm bau}_t)=0$. In other words, the emission market leads to a reduction of the production, and therefore a reduction of the carbon emission.

Let us summarize the present context of a small firm: (1)  the trading activity of the firm has no impact on its optimal production policy $q^{(1)}$ which is obtain from maximizing the profit of the firm,  (2) the firm's optimal production $q^{(1)}$ is smaller than that of the business-as-usual situation, so that the emission market is indeed a good tool for the reduction of carbon emission, and (3) the emission market assigns a price to the externality that the firm manager can use in order to optimize his production scheme.

\begin{Remark}\label{taxation}
Let us examine the case where there is no possibility to trade the carbon emission allowances, i.e. a standard taxation system where $\alpha$ is the amount of tax to be paid at the end of period per unit of carbon emission. Assuming again that the firm's horizon coincides with this end of period, its objective is:
 \b*
 V_0
 &:=&
 \sup_{q_.\in\Qc} \E\Bigl[ U\Bigl(\int_0^T \pi_t(q_t)dt-\alpha\left(E^q_T-E^{(0)}\right)_+\Bigr)\Bigr]
 \e*
Direct calculation leads to the following characterization of the optimal production level:
 \begin{equation}\label{q0}
 \begin{split}
 \frac{\partial\pi_t}{\partial q}\Bigl(q^{(0)}_t\Bigr)
 &=
 \alpha\frac{\partial \eta_t}{\partial q}\Bigl(q^{(0)}_t\Bigr)
 \E^{\Q^{(0)}}_t\Bigl[\Bigl(E_t^{q^{(0)}}-E^{(0)}\Bigr)_+\Bigr]\\
 \text{ with }\;\frac{d\Q^{(0)}}{d\P}
 &\propto
U'\Bigl(\int_0^T \pi_t(q^{(0)}_t)dt-\alpha\Bigl(E_t^{q^{(0)}}-E^{(0)}\Bigr)_+\Bigr).
  \end{split}
  \end{equation}
The natural interpretation of \eqref{q0} is that the production firm assigns an individual price to its emission:
 \begin{equation}\label{eqn:subjective_price}
 S_t :=
 \alpha\E^{\Q^{(0)}}_t\left[\left(E_t^{q^{(0)}}-E^{(0)}\right)_+\right],
 \end{equation}
i.e. the expected value of the amount of tax to be paid under the measure $\Q^{(0)}$ defined by its marginal utility as a density. The probability measure $\Q^{(0)}$ is the objective risk-neutral measure of the firm. Given this evaluation, the firm optimizes its adjusted profit function, $\pi_t(q)-\eta_t(q)S_t$; i.e. $\frac{\partial \pi_t}{\partial q}(q^{(0)})=\frac{\partial \eta_t}{\partial q}(q^{(0)})S_t$.
We continue by commenting on the optimal production policy defined by \eqref{q0}:
\begin{itemize}[leftmargin=*]
\item This problem would be considerably simplified if the manager were to know the market price for carbon emission. But of course, in the present context, \eqref{eqn:subjective_price} gives the firm's subjective price which is not quoted on any financial market and is hard to evaluate as the system of equations \eqref{q0} is still a nontrivial nonlinear fixed point problem.
\item The present situation, based on a classical taxation policy, offers no incentive to reduce emission beyond $E^{(0)}$. Indeed, if the optimal production is already below the level $E^{(0)}$, then it is indeed the same as the business-as-usual situation. So, the taxation does not contribute to further reduce the carbon emission.
\end{itemize}
\end{Remark}

The emission market provides an evaluation of the externality of carbon emission by firms. Given this information there is no more need to know precisely the utility function of the firm in order to solve the nonlinear system \eqref{q0}. The quoted price of the externality is then very valuable for the managers as it allows them to better optimize their production scheme.

%%%%%%%%%%%%%%%%%%%%%%%%%%%%%%%%%%%%%%%%%%%%%%%%%%%%
%%%%%%%%%%%%%%%%%%%%%%%%%%%%%%%%%%%%%%%%%%%%%%%%%%%%
%%%%%%%%%%%%%%%%%%%%%%%%%%%%%%%%%%%%%%%%%%%%%%%%%%%%
\section{Large producer with one-period carbon emission market}
\label{sec:large}
In this section, we consider the case of a large carbon emitting production firm. We shall see that this leads to different considerations as the trading activity  have an impact on the production policy of the firm. 
We model this situation by assuming that the state variable $Y$ is affected by the production policy of the firm:
$ dY_t^q
 =
 \left(\mu_t+\beta \eta_t(q_t)\right)dt + \gamma_t dW_t
$
where $\beta>0$ is a given impact coefficient. The price process $S$ of the carbon emission allowances is, as in the previous section, given by the no-arbitrage valuation principle $S_t^q \;=\; \alpha\Q_t^q\left[Y_T^q\ge 0 \right]$ and is also affected by the production policy $q$. The equivalent martingale measure $\Q^q$ is characterized by its Radon-Nykodim density which can be represented as a Dol\'eans-Dade exponential martingale generated by some risk premium process $\lambda$.
\[ 
\left.\frac{d\Q^q}{d\P}\right|_{\Fc_t}= \exp\Bigl( -\int_0^t \lambda_s(q_s)dW_s-\frac12\int_0^t \lambda_s(q_s)^2ds\Bigr)
\] 
The martingale property of of the bove Dol\'eans-Dade exponential follows from Assumption \ref{assp:pi_eta_lambda}-(iv) presented later.
In general, the risk premium process $\lambda$ may depend on the path of the control process $q$. For technical reasons, we shall restrict the analysis to those risk-neutral probability measures with risk premium process depending on the current value of the control process. i.e. $\lambda:[0,T]\times\Omega\times\R_+\longrightarrow\R$ is an $\F-$progressively measurable map. Under $\P$, the dynamics of the price process $S$ is given by $\frac{dS_t^q}{S_t^q}= \sigma^q_t\left(dW_t+\lambda_t(q_t)dt\right)$ for $t<T$,
where the volatility function $\sigma^q_t$ is progressively measurable and depends on the control process $\{q_s,0\le s\le T\}$. 
\begin{Remark}\label{rem:permit price}
The study performed in \cite{smw08} supports the assumption of existence of a martingale measure. In fact, by using empirical data, they showed that  the discounted price of the allowance is martingale. As a consequence, there is no seasonal effect in the price and we can simply assume that $\sigma^q$ is independent of time $t$. 
\end{Remark}

The effect of large producer on the market price of allowances is two-fold; one by directly adding to the drift of process $Y^q$ and the other by impacting the way the market evaluates the allowances, i.e. by changing martingale measure $\Q^q$. To separate the analysis of these two effects, in the next section we first consider the case where the risk premium of the market is not affected by the large producer. 

\subsection{Large carbon emission with no impact on risk premium}\label{sec:large_no_impact}

In this subsection, we restrict our attention to the case of large emitting firm with no impact on the risk premium, i.e. $
 \lambda_t(q)$ is independent of $q$ for $t\ge 0$.
The objective of the large emitting firm is:
 \b*
 V^{(2)}_0
 &:=&
 \sup_{q_\cdot\in\Qc,~\theta\in\Ac} \E\left[ U\left(X^{x,\theta}_T + B^q_T \right) \right].
 \e*

\begin{Proposition}\label{proplargeemission1}
Assume that the risk premium $\lambda$ is independent of $q$. Then, the optimal production policy is independent of the utility function of the producer $U$, and obtained by solving 
 \begin{equation}\label{prob:max_q2}
\sup_{q_\cdot\in\Qc} \E^\Q\left[ B^q_T \right]
 \end{equation}
 where $\Q$ is the martingale measure.
Moreover, if $q^{(2)}$ is an optimal production scheme, then the optimal investment strategy $\theta^{(2)}$ is characterized by
 \begin{equation}\label{eqn:optimal_theta_large}
 X_T^{x,\theta^{(2)}}+B_T^{q^{(2)}}=(U')^{-1}\left(y^{(2)}\frac{d\Q}{d\P}\right),
 ~~
 x+\E^\Q\left[B_T^{q^{(2)}}\right]
 =
 \E^\Q\left[(U')^{-1}\left(y^{(2)}\frac{d\Q}{d\P}\right)\right].
 \end{equation}
\end{Proposition}
\begin{proof} The proof follows the same line of arguments as in Proposition \ref{prop:small_theta_q}.
\end{proof}      

\vspace{5mm}

In order to push further the characterization of the optimal production policy $q^{(2)}$, we specialize the discussion to the Markov case by assuming the following for the triple $(\pi,\eta,\lambda)$.
\begin{assp}{A}\label{assp:pi_eta_lambda}
$\pi_t(q)=\pi_t(q)$, $\eta_t(q)=\eta_t(q)$, and $\lambda_t(q)=\lambda(t)$ are in $\Cn{0,1}([0,T]\times\R_+)$ and satisfy
\begin{enumerate}[label={(\roman*) }, leftmargin=*]
\item $\pi$ is strictly concave in $q$, $\pi_t(0)=0$ and $\frac{\partial\pi}{\partial q}(\infty)<0$, 
\item $\eta$ is convex and strictly increasing in $q$,
\item $\lambda$ is concave and nondecreasing in $q$ and $\lambda(t,0)\ge0$,
\item $|\sup_{q\ge0}\{\pi+\eta_t(q)v_2 +(\eta_t(q)-\gamma_t(y)\lambda_t(q))v_1\}|\le C|v_1|^2+g(t,y)$ for some $C>0$.% and {\Large\color{red} Check Wang II}.
\end{enumerate}
\end{assp} 
We also enforce a Markovian dynamics for process $Y^q$ under measure $\P$; i.e.
\b*
dY_t^q=\left(\mu_t(Y_t^q)+\beta \eta_t(q_t)\right)dt + \gamma_t(Y_t^q) dW_t,
\e*
for some deterministic functions $\mu,\gamma:[0,T]\times\R\longrightarrow\R$ and a nonnegative constant $\beta$.

The controlled variable $E_t^q$ is now defined by the dynamics $dE_t^q = \eta_t(q_t)dt$
which records the cumulated carbon emission of the firm. The dynamic version of the production policy optimization problem \eqref{prob:max_q2} is given by:
\[
 V^{(2)}(t,y,e)
 :=
 \sup_{q_\cdot\in\Qc_t} \E^\Q_{t,e,y} \left[ \int_t^T \pi(s,q_s)ds
                                      -\alpha{1}_{\{Y^q_T\ge0\}}E^q_T \right],
 \]
  where $\Qc_t$ is the collection of all non-negative $\F-$progressively measurable processes such that $\int_t^T\eta_s(q_s)ds<\infty$, and $\E^\Q_{t,e,y}$ is the  expectation with respect to $\Q$ conditional on $E^q_t=e$, $Y^q_t=y$. Here we absorb the initial free allowances $E^{(0)}$ into the  condition $E^q_t=e$ by assuming that $e$ can take negative values.
Then, $V^{(2)}$ is a viscosity solution of the dynamic programming equation with a terminal condition :
 \begin{equation} \label{eqn2}
 \begin{split}
 0
 &=
 -\frac{\partial V^{(2)}}{\partial t}-(\mu-\lambda\gamma) V^{(2)}_y - \frac12 \gamma^2 V^{(2)}_{yy}-\sup_{q\ge 0} \theta(t,q,V^{(2)}_y,V^{(2)}_e)\\
& V^{(2)}(T,y,e)= -\alpha{1}_{\{y>0\}}e,
 \end{split}
 \end{equation}
 where 
$\theta(t,q,v_1,v_2)= \pi_t(q)+\eta_t(q)v_2+\beta \eta_t(q)v_1$.
By Lemmas \ref{lem:regularity_in_e} and \ref{lem:regularity_in_t_y} and Corollary \ref{cor:right_left_derivatives}, $V^{(2)}$ is continuously differentiable once in $t$ and twice in $y$ when $t<T$, Lipschitz in $e$ and $\partial_{e+}V$ exists and is right-continuous. Then, the optimal production $q^{(2)}$ is given by the maximum 
\begin{equation}\label{eqn:characterization_of_q^(2)}
q^{(2)}(t,y,e)\in\argmax_{q\ge 0} \left\{ \pi_t(q)+\eta_t(q)(V^{(2)}_{e+}+\beta V^{(2)}_y)(t,y,e) \right\}
\end{equation}
By Lemma \ref{lem:V_e+=-S_t}, we have
\be\label{eqn:V^2_e+=-S_t}
-V^{(2)}_{e+}(t,\eta_t,Y_t)=S_t.
\ee
If the maximum in \eqref{eqn:characterization_of_q^(2)} is attained in an interior point, then it satisfies
 \be\label{comparison2}
 \frac{\partial \pi}{\partial q}\left(t,q^{(2)}_t\right)
 &=&
 \frac{\partial \eta}{\partial q}\left(t,q^{(2)}_t\right) 
 \left(S_t-\beta V^{(2)}_y(t,E^{q^{(2)}}_t,Y^{q^{(2)}}_t)\right).
 \ee
Otherwise in case the maximum in \eqref{eqn:characterization_of_q^(2)} is not attained in an interior point, we have $q^{(2)}=0$.
Thus, it follows from comparing \eqref{comparison2} with $ q^{(2)} \le q^{(1)}$ if and only if $V^{(2)}_y\le0$. In fact, a larger positive $-\beta V^{(2)}_y(t,E^{q^{(2)}}_t,Y^{q^{(2)}}_t)$ implies a smaller $q^{(2)}$ below $q^{(1)}$. For instance if $E_0^{q^{(2)}}\ge0$, $V^{(2)}_y$ remains non-positive at all time. In this case, by choosing a large penalty term $\alpha$, the optimal production and consequently the emission can be controlled to meet the target.
In other words, the impact of the production of the firm on the prices of carbon emission allowances increases the cost of the externality for the firm. This immediately affects the profit function of the firm and leads to a decrease of the level of optimal production. Hence, the presence of the emission market is playing a positive role in terms of reducing the carbon emission. 
The following result summarizes the above discussion.
\begin{Theorem}\label{thm:largenoimpact}
Let Assumption \ref{assp:regularity_Y} holds and triple $(\pi,\eta,\lambda)$ satisfy Assumption \ref{assp:pi_eta_lambda}.
 Then, $V^{(2)}$ satisfies problem \eqref{eqn2}  and is in $\Cn{1,2,0}([0,T)\times\R\times\R)$, $\partial_{e+}V^{(2)}$ exists and is right continuous. \\
In addition, \ref{eqn:V^2_e+=-S_t} holds and optimal production policy is characterized by \eqref{eqn:characterization_of_q^(2)}.
\end{Theorem}
\begin{Corollary}\label{cor:largenoimpact}
Under the same assumption as Theorem \ref{thm:largenoimpact}, $q^{(2)}\le q^{(1)}$.
\end{Corollary}

\subsection{Large carbon emission impacting the risk-neutral measure}

We now consider the general case where the risk premium process is impacted by the emission of the production firm:
 \b*
 \left.\frac{d\Q^q}{d\P}\right|_{\Fc_T}
 &=&
 \exp\left(-\int_0^T \lambda(q_t)dW_t-\frac12\int_0^T \lambda(q_t)^2dt \right).
 \e*
The partial maximization with respect to $\theta$, as in the proof of Proposition \ref{proplargeemission1}, is still valid in this context, and reduces the production firm's problem to
 \be\label{prob:large2}
 \sup_{q_\cdot\in\Qc} \E\left[ U\circ(U')^{-1}\left(y^q\;\frac{d\Q^q}{d\P}\right)\right]
 \ee
where $y^q$ is defined by 
 \be\label{budget_constraint}
 \E^{\Q^q}\left[(U')^{-1}\left(y^q\;\frac{d\Q^q}{d\P}\right)\right]
 &=&
 x+\E^{\Q^q}\left[B^q_T\right].
 \ee 
In order to move further, we assume that the preferences of the production firm are defined by an exponential utility function
 \b*
 U(x) &:=& -e^{-a x},~~x\in\R~~a>0.
 \e* 
Then $U\circ(U')^{-1}(y)=-y/a$, and \eqref{prob:large2} reduces to 
 \be\label{prob:large3}
 \inf_{q_.\ge 0}\E\left[y^q\;\frac{d\Q^q}{d\P}\right]
 &=&
 \inf_{q_.\ge 0}\;y^q.
 \ee
Finally, the budget constraint \eqref{budget_constraint} is in the present case:
 \b*
 x+\E^{\Q^q}\left[B^q_T\right]
 &=&
 \frac{-1}{a}\E^{\Q^q}\left[\ln{\left(\frac{y^q}{a}\;\frac{d\Q^q}{d\P}\right)}\right]
 \\
 &=&
 \frac{-1}{a}\left\{\ln{\left(\frac{y^q}{a}\right)}
                      +\E^{\Q^q}\left[\ln{\left(\frac{d\Q^q}{d\P}\right)}\right]
                \right\},
 \e*
so that the optimization problem \eqref{prob:large3} is equivalent to:
 \begin{equation}
\sup_{q_\cdot\in\Qc} \E^{\Q^q}\left[B^q_T +\frac{1}{a}\ln{\left(\frac{d\Q^q}{d\P}\right)} \right]
 =
 \sup_{q_\cdot\in\Qc} \E^{\Q^q}\left[\int_0^T \left(\pi+\frac{\lambda^2}{2a}\right)(t,q_t)dt-S^q_TE_T^q\right].
 \label{prob:large4}
 \end{equation}
 Notice the difference between the above optimization problem, which determines the optimal production policy of the production firm, and the problem in Section \ref{sec:large_no_impact} where the firm does not impact the risk premium. In the present section, the firm's optimization criterion is penalized by the entropy of the risk-neutral measure with respect to the statistical measure. Unlike Section \ref{sec:large_no_impact}, the optimal production of the firm with impact on the risk premium of the market depends on the risk preference of the firm.

The firm's optimal production problem \eqref{prob:large4} is a standard stochastic control problem. We continue this discussion by 
considering the Markov case, and introducing the dynamic version of \eqref{prob:large4}:
 \begin{equation}\label{eqn:large_lambda_impact}
 V^{(3)}(t,y,e)
 :=
 \sup_{q_\cdot\in\Qc_t} \E^{\Q^q}_{t,e,y}\left[\int_t^T\left(\pi+\frac{\lambda^2}{2a}\right)(s,q_s)ds
                                -\alpha{1}_{\{Y_T^q\ge 0\}E^q_T}
                          \right],
 \end{equation} 
where the controlled state dynamics is given by:
 \begin{equation}\label{eqn:Y_under_Q^q}
 \begin{split}
 dY_t^q &= \left(\mu_t(Y_t^q)+\beta \eta_t(q_t)
                -\gamma_t(Y_t^q)\lambda_t(q_t)\right)dt + \gamma_t(Y_t^q) dW^q_t,
 \\
 dE_t^q &= \eta_t(q_t)dt,
 \end{split}
 \end{equation}
$W^q$ is a Brownian motion under $\Q^q$, and $\mu$ and $\gamma$ are $\Cn{1,2}$ functions in $(t,y)$, and $\mu$, $\eta$ and $\lambda$ are $\Cn{1,2}$ functions in $(t,q)$.

By classical arguments, we then see that $V^{(3)}$ is a viscosity solution of
 \begin{equation}\label{eqn:V3}
 \begin{cases}
 0
 =
 \partial_t V^{(3)}+\mu V^{(3)}_y + \frac12\gamma^2 V^{(3)}_{yy} 
 +\mathop{\max}\limits_{q\in\R_+} \theta(t,y,q,V^{(3)}_e,V^{(3)}_y)&\text{\rm on}\;[0,T)\times\R^2\\
 V^{(3)}(T,y,e)
 =
 -\alpha {1}_{\{y>0\}} e&\text{\rm on }\;\R^2,
 \end{cases}
 \end{equation}
 where 
 \[
 \theta(t,y,q,p_e,p_y)= \pi_t(q)+\frac{1}{2a}\lambda_t(q)^2+\eta_t(q)(p_e+\beta p_y)-\gamma_t(y)\lambda_t(q)p_y
 \]
In terms of the value function $V^{(3)}$, the optimal production policy is obtained as the maximizer in the above equation, i.e.
\begin{equation}\label{eqn:characterization_of_q^(3)}
\begin{split}
q^{(3)}(t,y,e)&\in\argmax_{q\ge 0} \biggl\{ \pi_t(q)+\frac{1}{2a}\lambda_t(q)^2+\eta_t(q)(V^{(3)}_e+\beta V^{(3)}_y)(t,y,e)\\
&\hspace*{6cm}-\gamma_t(y)\lambda_t(q)V^{(3)}_y(t,y,e) \biggr\}.
\end{split}
\end{equation} 
Observe that if we assume $V^{(3)}$ is regular enough, then Assumption \ref{assp:pi_eta_lambda} implies that argmax is a singleton and $q^{(3)}$ is unique.
In addition if an interior maximum occurs, then the first order condition is:
\[ \frac{\partial \pi}{\partial q}(q^{(3)})
 +\frac{1}{a}(\lambda\frac{\partial \lambda}{\partial q})(q^{(3)})
 +\frac{\partial \eta}{\partial q}(q^{(3)})(V^{(3)}_e+\beta V^{(3)}_y)
 -\gamma\frac{\partial \lambda}{\partial q}(q^{(3)})V^{(3)}_y
 =
 0.
\]
Moreover, we shall show in Lemma \ref{lem:V_e+=-S_t} that the price of the carbon emission allowance, as observed on the emission market, is given by:
 \begin{equation}\label{eqn:V^3_e+=-S_t}
 S_t
=
 -V_{e+}^{(3)}(t,E^{q^{(3)}}_t,Y^{q^{(3)}}_t).
\end{equation}
 Then, it follows that the optimal production policy of the firm is defined by:
 \be
 \frac{\partial \pi}{\partial q}(t,q^{(3)})
 &=&
 \frac{\partial \eta}{\partial q}(t,q^{(3)})\left(S_t-\beta V_y^{(3)}(t,E^{q^{(3)}}_t,Y^{q^{(3)}}_t)\right)
 \nonumber
 \\
 &&
 +\frac{\partial \lambda}{\partial q}(t,q^{(3)})
   \left(\gamma V^{(3)}_y(t,E^{q^{(3)}}_t,Y^{q^{(3)}}_t)-\frac{1}{a}\lambda_t(q^{(3)})\right).
 \label{comparison3}
 \ee
Contrary to the previous case where the risk premium process was not impacted by the carbon emission of the large firm, we cannot always conclude from the above formula that $q^{(3)}$ is smaller than $q^{(1)}$; the optimal production policy in the absence of a financial market given by \eqref{defn:q^1}.
More precisely, if $tau$ defined below is non-negative, then we can conclude that $q^{(3)}\le q^{(1)}$.
 \begin{equation}\label{eqn:tau}
 \begin{split}
 \tau&:=\left(\beta \frac{\partial \eta}{\partial q}(t,q^{(3)}) -\gamma\frac{\partial \lambda}{\partial q}(t,q^{(3)})\right) V^{(3)}_y(t,E^{q^{(3)}}_t,Y^{q^{(3)}}_t)+\frac{1}{a}\frac{\partial \lambda}{\partial q}(t,q^{(3)})\lambda_t(q^{(3)})
 \end{split}
 \end{equation}
However, $\tau$ has no known sign, and there is no economic argument supporting that it should have some specific sign. Under Assumption \ref{assp:pi_eta_lambda}, we can only be sure that $\frac{1}{a}\frac{\partial \lambda}{\partial q}(t,q^{(3)})\lambda_t(q^{(3)})\ge0$. However, while $V^{(3)}_y\le0$, $\beta \frac{\partial \eta}{\partial q}(t,q^{(3)}) -\gamma\frac{\partial \lambda}{\partial q}(t,q^{(3)})$ does not have a known sign. Therefore due to the impact on the emission market, the optimal production of the  large producer can potentially be higher than the case when there is no emission market. Based on the discussion above, the case where we can make sure $q^{(3)}\le q^{(1)}$ is provided in the following proposition.
The above discussion is made rigorous in the following results which follows from Appendix \ref{sec:analytics}.
\begin{Theorem}\label{thm:largeimpact}
Let Assumption \ref{assp:regularity_Y} holds and triple $(\pi+\frac{\lambda^2}{2},\eta,\lambda)$ satisfy Assumption \ref{assp:pi_eta_lambda}.
 Then, $V^{(3)}$ is a $\Cn{1,2,0}([0,T)\times\R\times\R)$ solution of  problem \eqref{eqn:V3},  $\partial_{e+}V^{(3)}$ exists and is right continuous. \\
In addition, \ref{eqn:V^3_e+=-S_t} holds and optimal production policy is characterized by \eqref{eqn:characterization_of_q^(3)}.
\end{Theorem}
\begin{Corollary}\label{cor:largeimpact}
Under the same assumption as Theorem \ref{thm:largenoimpact}, if we have 
$\beta \frac{\partial \eta}{\partial q}(t,q^{(3)}) -\gamma\frac{\partial \lambda}{\partial q}(t,q^{(3)})\le0$, then $q^{(3)}\ge q^{(1)}$.
\end{Corollary}
   
In the next section, we discuss the cases where $\beta \frac{\partial \eta}{\partial q}(t,q^{(3)}) -\gamma\frac{\partial \lambda}{\partial q}(t,q^{(3)})>0$ through numerical implementation of HJB equation \eqref{eqn:V3} to determine the region where $q^{(3)}>q^{(1)}$. An important question in this case is how much the total emission is affected for different choice of parameters $\gamma$ and $\alpha$ controlled by the regulator.
%%%%%%%%%%%%%%%%%%%%%%%%%%%%%%%%%%%%%%%%%%%%%%%%%%%%
%%%%%%%%%%%%%%%%%%%%%%%%%%%%%%%%%%%%%%%%%%%%%%%%%%%%
%%%%%%%%%%%%%%%%%%%%%%%%%%%%%%%%%%%%%%%%%%%%%%%%%%%%
\begin{comment}
\section{Multiperiod model with banking}
The analysis of the previous sections are restricted to the case where the carbon allowances market is organized over one single period. In this section, we discuss how to extend this results to a multiperiod model with {banking}. We then assume that there are $n$ maturities for the carbon allowances market
 $$
 T_1<\ldots <T_n
 $$
instead of a single one. According to the banking rule, the carbon emission allowance can serve for the next periods if not used for the current one. Then, the allowance can be viewed as a derivative security with payoff:
 $$
 S_{T_n}
 :=
 \alpha\left({1}_{\{Y_{T_1}^q\ge0\}}
             +{1}_{\{Y_{T_1}^q<0\}}{1}_{\{Y_{T_2}^q\ge0\}}
             +\ldots
             +{1}_{\{Y_{T_i}^q<0,~i<n\}}
              {1}_{\{Y_{T_n}^q\ge0\}}
        \right).
 $$
Following the same argument as in the previous section, the no-arbitrage market price at each time $t\le T_n$ is given by 
 \be\label{S-banking}
 S_t
 :=
 \E_t^\Q\left[S_{T_n}\right]
 &\mbox{for all}&
 t\le T_n,
 \ee
where $\Q$ is the risk-neutral measure. Now it is clear that all the analysis of the previous sections apply by just replacing the price formula \eqref{S-1period} by the above market price \eqref{S-banking}.
\end{comment}
%%%%%%%%%%%%%%%%%%%%%%%%%%%%%%%%%%%%%%%%%%%%%%%%%%%%
%%%%%%%%%%%%%%%%%%%%%%%%%%%%%%%%%%%%%%%%%%%%%%%%%%%%
%%%%%%%%%%%%%%%%%%%%%%%%%%%%%%%%%%%%%%%%%%%%%%%%%%%%
\section{Numerical results}
\label{Sec:numerics}
The main goal of the numerical results  is to understand the behavior of the optimal strategy $q^{(3)}$ in \eqref{comparison3} and more precisely to study the case where $q^{(3)}>q^{(1)}$.
If we consider $\pi(q)=q(1-q)$, $\eta(q)=\lambda(q)=q$, $\beta=1$, $\gamma=.65$, $T=10$ and at this moment $\alpha=0.1$, 
then \eqref{eqn:V3} reduces to
\begin{equation}\label{eqn:pde_example}
\partial_tV+\mu V_y+\frac{1}{2}\gamma^2V_{yy}+\frac{1}{4\varrho}\left(1+V_{e}+(1-\gamma)V_y \right)_+^2=0. 
\end{equation}
Note that %this example satisfies Assumption \ref{assp:pi_eta_lambda} except (iv). However by Remark \ref{rem:regularity_linear}, we can still have sufficient regularity to conclude 
by Lemmas \ref{lem:optimal_control} and \ref{lem:V_e+=-S_t}, we have 
$V_{e+}=-S_t$ and optimal control is given by
\[
q^{(3)}=\frac{1}{2\varrho}(1+V_e+(1-\gamma)V_y)_+,
\]
where $\varrho=\left(1-\frac{1}{2a}\right)$, and we used direct calculations to obtain
\[
\mathop{\max}\limits_{q\geq 0}\theta(t,y,q,V_e,V_y)=\frac{1}{4\varrho}\left(1+V_e+(1-\gamma)V_y \right)_+^2.
\]
To determine the region where $q^{(3)}>q^{(1)}$, we have to find the region where $\tau$  given by \eqref{eqn:tau} is positive, i.e.
\[
 \tau
 := 
(1-\gamma)V_y+(\varrho^{-1}-1)(1+V_{e+}+(1-\gamma)V_y)_+>0.
 \]
For the choice of parameters 
$\mu=0.1$, $\rho=.9$ ($a=5$) and $T=10$, we approximated the value function, correction term $\tau$, and optimal control by a finite-difference Trotter-Kato based scheme whose details is given in Appendix \ref{appnedix_numerics}.
\begin{figure}
\begin{subfigure}[b]{0.3\textwidth}
\centering
\resizebox{\linewidth}{!}{
\includegraphics[scale=1]{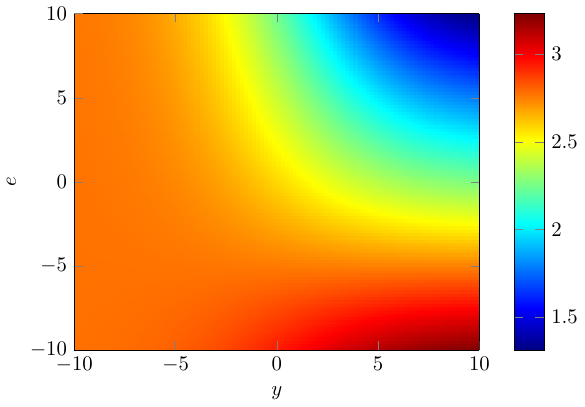};
}
\caption{Solution $V(0,\cdot)$}
\end{subfigure}
\begin{subfigure}[b]{0.32\textwidth}
\centering
\resizebox{\linewidth}{!}{
\includegraphics[scale=1]{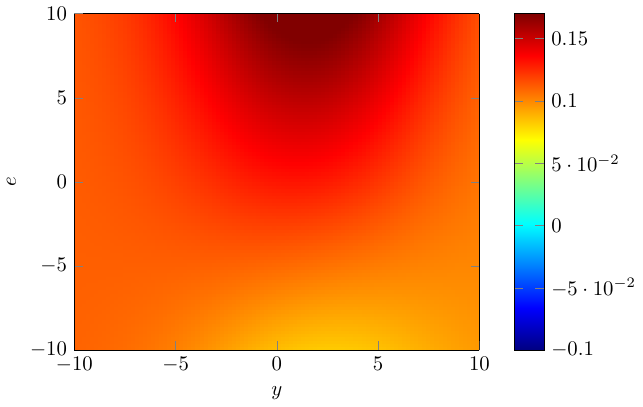};
}
\caption{Correction term $\tau$}
\end{subfigure}
\begin{subfigure}[b]{0.3\textwidth}
\centering
\resizebox{\linewidth}{!}{
\includegraphics[scale=1]{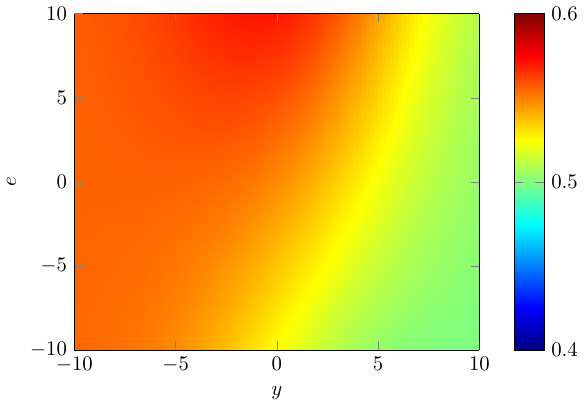};
}
\caption{Optimal Control $q^{(3)}$}
\end{subfigure}
\caption{When $\gamma=1.5$ and $\alpha=0.1$}
\label{fig:gamma=1.5}
\end{figure}

%%%%%%%%%%%%%%%%%%%%%%%%%%%%%%
%%%%%%%%%%%%%%%%%%%%%%%%%%%%%%

\begin{figure}
\begin{subfigure}[b]{0.3\textwidth}
\centering
\resizebox{\linewidth}{!}{
\includegraphics[scale=1]{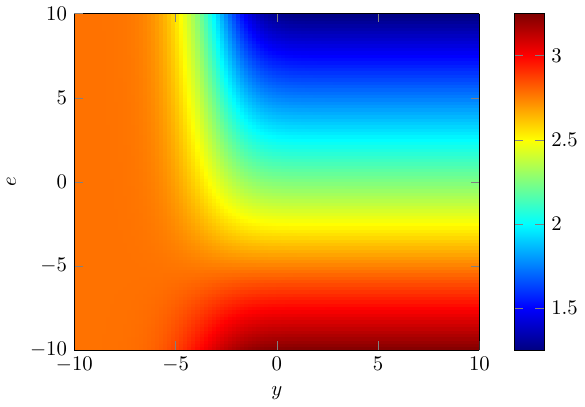};
}
\caption{Solution $V(0,\cdot)$}
\end{subfigure}
\begin{subfigure}[b]{0.32\textwidth}
\centering
\resizebox{\linewidth}{!}{
\includegraphics[scale=1]{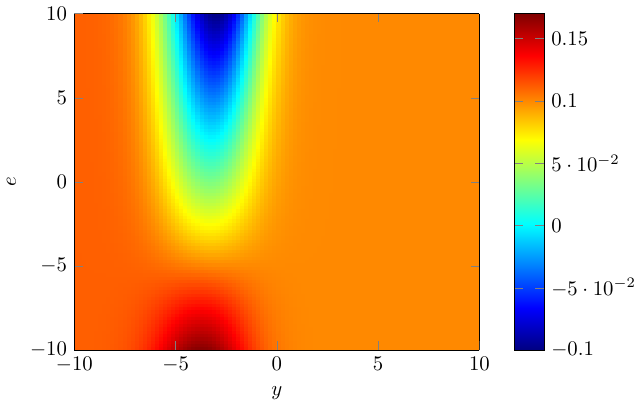};
}
\caption{Correction term $\tau$}
\end{subfigure}
\begin{subfigure}[b]{0.3\textwidth}
\centering
\resizebox{\linewidth}{!}{
\includegraphics[scale=1]{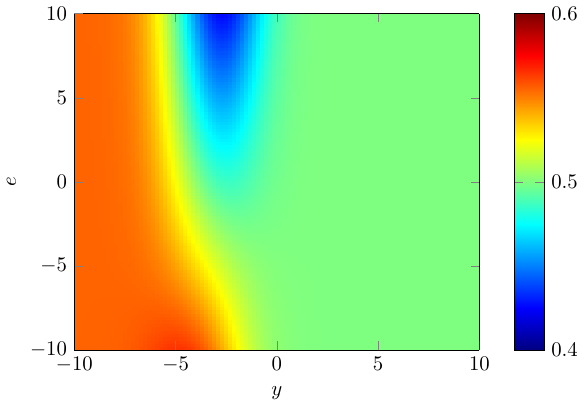};
}
\caption{Optimal Control $q^{(3)}$}
\end{subfigure}
\caption{When $\gamma=0.5$ and $\alpha=0.1$}
\label{fig:gamma=0.5}
\end{figure}

As shown in Proposition \ref{cor:largeimpact} and Figure \ref{fig:gamma=1.5}, the optimal production $q^{(3)}$ for large production firm with impact on risk premium is always higher that $q^{(1)}$ when $\gamma=1.5$.  However when $\gamma=0.5$, there is a region where we have $q^{(3)}\le q^{(1)}$ which is shown in Figure  \ref{fig:gamma=0.5} for $t=0$. 

%As a matter of fact, the position of the large producer can be distanced from region $\{\tau>0\}$ by having enough of allowance in the portfolio at time $0$. On the other hand, too much uncertainty about the process $Y$ measure with the value of $\gamma$ can eventually put the producer in distress, since the size of this region grows as $\gamma$ becomes larger. Therefore, we can avoid large producer from producing more than $q^{(1)}$ by giving sufficient free allowances or/and reducing the uncertainty of the pollution index by frequent  announcements.

%%%%%%%%%%%%%%%%%%%%%%%%%%%%%%%%%%%%%%%%%%%%%%%%%%%%
%%%%%%%%%%%%%%%%%%%%%%%%%%%%%%%%%%%%%%%%%%%%%%%%%%%%
%%%%%%%%%%%%%%%%%%%%%%%%%%%%%%%%%%%%%%%%%%%%%%%%%%%%

\appendix

%%%%%%%%%%%%%%%%%%%%%%%%%%%%%%%%%%%%%%%%%%%%%%%%%%%%
%%%%%%%%%%%%%%%%%%%%%%%%%%%%%%%%%%%%%%%%%%%%%%%%%%%%
%%%%%%%%%%%%%%%%%%%%%%%%%%%%%%%%%%%%%%%%%%%%%%%%%%%%
\section{Uniqueness, verification and existence of optimal control}
\label{sec:analytics}
Throughout the appendix, we assume that $(\Omega,\Fc,\F=\{\Fc_t\}_{t\ge0},\P)$ is a filtered probability space satisfying the usual conditions which hosts a Brownian motion $\{W_t\}_{t\ge0}$ and let $\E$ denote the expectation with respect to $\P$. Let
\begin{equation}\label{prob:control_problem}
\begin{split}
V(t,y,e)
&=
\sup_{q_\cdot\in\Qc_t}J_{q}(t,y,e)\\
J_{q}(t,y,e)&=\E_{t,e,y}\left [\int_t^T\tilde\pi(s,q_s)\diff s-\alpha{1}_{\{Y_T^{q,y}\ge0\}}E_T^{q,e}\right ],
\end{split}
\end{equation}
where $dE_t^q=\eta_t(q_t)dt$ and
\begin{equation}\label{eqn:Y^q}
dY_t^q=\Bigl(\mu_t(Y_t^q)+\beta \eta_t(q_t)-\gamma_t(Y_t^q)\lambda_t(q_t)\Bigr)dt+\gamma_t(Y_t^q)dW_t,\;
\end{equation}
where  
$\mu,\gamma:[0,T]\times\R\to\R$ are continuous in $t$ and  Lipschitz in $y$ with $\gamma\ge0$.
\begin{Remark}\label{rem:Q super q in dynamics of Y}
Notice that  in the current Appendix, the reference probability  measure $\P$ is  different from the physical probability measure introduced at the beginning of Section \ref{Sec:small}. This setting helps us extend the results in this appendix to both value functions $V^{(2)}$ and $V^{(3)}$. More precisely, if we set $\tilde\pi=\pi$ and  $\P=\Q$, then $V=V^{(2)}$. Else if $\tilde \pi=\pi+\frac{\lambda^2}{2a}$ and $\P=\Q^q$, then $V=V^{(3)}$; here the dependency of martingale measure  $\Q^q$ with respect to $q$ in the definition of $V^{(3)}$ is absorbed in the dynamic of $Y_t^{q}$. 
\end{Remark}
We would like to show that $V$ can be characterized by the HJB equation
 \begin{equation} \label{prob:HJB_Appendix}
 \begin{cases}
 0
 =
-\partial_t V-H(t,y,\partial_y V,\partial_{e} V, \partial_{yy} V)&\;\text{\rm for}\;(t,e,y)\in[0,T)\times\R\times\R\\
 V(T,y,e)
 =
 -\alpha{1}_{\{y>0\}} e& \;\text{\rm for}\;(e,y)\in\R\times\R,
 \end{cases}
 \end{equation}
where
\begin{equation*}
\begin{split}
H(t,y,v_1,v_2,v_{11})&:=\mu_t(y)v_1 + \frac12\gamma_t^2(y)v_{11}+\sup_{0\le q} \theta(t,y,v_1,v_2,v_{11}) \\
 \theta(t,y,v_1,v_2,v_{11})&:=\tilde\pi_t(q)+\eta_t(q)v_2+(\beta\eta_t(q)-\gamma_t(y)\lambda_t(q))v_1 .
 \end{split}
\end{equation*}
Because of discontinuity in terminal condition, we adopt the definition of discontinuous viscosity solutions from \cite[Section 6.2]{Touzi12} or \cite[Section 4.2]{Barles94} . For a locally bounded measurable function $u$, we denote by $u^*$ and $u_*$ the upper semicontinuous and lower semicontinuous envelopes of $u$, respectively.
\begin{Definition}\label{defn:viscosity}
Let $h(y,e)$ be a locally bounded measurable function. For  \eqref{prob:HJB_Appendix} with terminal value $V(T,y,e)=h(y,e)$, 
\begin{enumerate}[label={(\alph*) }, leftmargin=*]
\item a locally bounded measurable  function $u$ upper semicontinuous on $[0,T)\times\R\times\R$ is called a viscosity subsolution  if
\begin{enumerate}[label={(\roman*) }, leftmargin=-.1cm]
\item $u^*(T,y,e)\le h^*(y,e)$
\item for any smooth function $\phi$ such that  $\max(u^*-\phi)=(u^*-\phi)(t_0,e_0,y_0)$ with $t_0<T$, at $(t_0,e_0,y_0)$ we have 
\begin{align*}
- \frac{\partial \phi}{\partial t}-H_*(t,y,\phi_y,\phi_e,\phi_{yy})\le 0
 \end{align*}
\end{enumerate}
\item a locally bounded measurable function $u$  lower semicontinuous on $[0,T)\times\R\times\R$ is called a viscosity supersolution of  \eqref{prob:HJB_Appendix} if
\begin{enumerate}[label={(\roman*) }, leftmargin=-.1cm]
\item $u_*(T,e',y')\ge h_*(y,e)$
\item for any smooth function $\phi$ such that  $\min(u_*-\phi)=(u_*-\phi)(t_0,e_0,y_0)$ with $t_0<T$, at $(t_0,e_0,y_0)$ we have 
\begin{align*}
-\frac{\partial \phi}{\partial t}-H^*(t,y,\phi_y,\phi_e,\phi_{yy})\ge 0
 \end{align*}
\end{enumerate}
\item a locally bounded measurable function $u$ continuous on $[0,T)\times\R\times\R$ is called a viscosity solution of  \eqref{prob:HJB_Appendix} if it is both a viscosity sub- and supersolution.
\end{enumerate}
\end{Definition}

 \begin{Theorem}\label{thm:viscosity}
Let triple $(\tilde\pi,\eta,\lambda)$ satisfy Assumption \ref{assp:pi_eta_lambda}. Then $V$ is the unique\footnote{The specific sense of uniqueness here is discussed in Remark \ref{rem: comparison at t=T}.} viscosity solution of \eqref{prob:HJB_Appendix} on $[0,T]\times\R\times\R$.
\end{Theorem}
\begin{proof} By Assumption \ref{assp:pi_eta_lambda}, $H$ is locally bounded, and can be approximated by a net of continuous functions $\{H_M\}_{M>0}$
\[
H_M(t,y,v_1,v_2,v_{11}):=\mu_t(y)v_1 + \frac12\gamma_t^2(y)v_{11} 
 +\sup_{0\le q\le M} \theta(t,y,v_1,v_2,v_{11}).
\] 
Thus, one can apply  \cite[Theorems 7.4 and 6.8]{Touzi12} to obtain (a.ii) and (b.ii) in Definition \ref{defn:viscosity}. To show (a.i) and (b.i), we approximate the terminal condition $h(y,e)=-\alpha e 1_{\{y\ge0\}}$ by two smooth  functions $-\alpha e\underline{\rho}_\eps(y)$ and $-\alpha e\overline{\rho}_\eps(y)$ from below and above respectively, i.e. $\underline{\rho}_\eps(y)=1$ on $y\ge 0$, $\underline{\rho}_\eps(y)=0$ on $y\le -\eps$ and $0\le \underline{\rho}_\eps(y)\le 1$, and $\overline{\rho}_\eps(y)=1-\underline{\rho}_\eps(-y)$.
Then by \cite[Theorems 7.4 and 7.6]{Touzi12}, the value functions  $\overline{V}_\eps$ and $ \underline{V}_\eps$ defined below are the unique continuous viscosity solutions\footnote{in class of functions with linear growth} of problem \eqref{prob:HJB_Appendix} with terminal condition $\underline{V}_\eps(T,y,e)= -\alpha e \underline{\rho}_\eps(y)$ and $\overline{V}_\eps(T,y,e)= -\alpha e \overline{\rho}_\eps(y)$, respectively.
 \begin{equation*}
\overline{ \underline{V}}_\eps(t,y,e)
 =
 \sup_{q_\cdot\in\Qc} \E_{t,e,y}\left[\int_t^T\tilde\pi_t(q_t)dt
                                -\alpha \underline{\overline{\rho}}_\eps(Y_T^q)E^{q}_T
                          \right]
 \end{equation*}  
 On the other hand,  it follows from the the optimal control problems of  $\underline{V}_\eps$, $V$, and $\overline{V}_\eps$ that 
 $\underline{V}_\eps\le V\le \overline{V}_\eps$. Therefore, by taking upper semicontinuous and lower semicontinuous envelopes from both sides and then sending $\eps\to 0$,  we obtain the desired result.\\
 For uniqueness, first  notice that $\overline{V}_\eps-\underline{V}_\eps\to0$ as $\eps\to0$. By standard comparison, e.g. \cite[Theorem 6.21]{Touzi12}, Since  any upper semicontinuous subsolution  $u$ (lower semicontinuous supersolution $v$) of \eqref{prob:HJB_Appendix}, is also an upper semicontinuous subsolution  (a lower semicontinuous supersolution)  of HJB problem for $\overline{V}_\eps$ ($\underline{V}_\eps$), we have $u\le \overline{V}_\eps$ ($v\ge \underline{V}_\eps$). Thus, 
 $u-v\le \overline{V}_\eps-\underline{V}_\eps$ and by sending $\eps\to0$, we obtain uniqueness for $t<T$.
\end{proof}      
\begin{Remark}\label{rem: comparison at t=T}
The above continuity result does not imply that $u(T,e,0)\le v(T,e,0)$. In fact, it only implies that $u\le v+\alpha e \delta_{0}(y)\delta_T(t)$. In this case, the uniqueness may be violated along the half-line $\{(T,e,0): e>0\}$. But, we can see that it does not affect the main results of this study.
\end{Remark}
Theorem  \ref{thm:viscosity} requires minimal regularity of the value function $V$. However to achieve the results of Section \ref{sec:large} for the large production firm, we need to show that an optimal control exists and can be expressed in terms of derivatives of $V$. To do so, we need to impose Assumption \ref{assp:regularity_Y}.
\begin{assp}{Y}\label{assp:regularity_Y}
$\mu,\gamma:[0,T]\times\R\longrightarrow\R$  are $\Cn{\infty,\infty}$ and there is some positive constant $c$, such that for all $(t,y)$ $\gamma_t(y)\ge c>0$.
\end{assp}
\begin{Remark}\label{rem:regulairty_in_(t,y)}The above assumption implies that the semigroup $\{P_t\}_{t\ge0}$, generated by operator $\Lc:=\frac{\gamma^2}{2}\partial_{yy}+\mu\partial_y$ on a bounded regular domain $Q\subseteq[0,\infty)\times\R$ with continuous boundary conditions, is in $\Cn{\infty,\infty}$, in the sense that for any  bounded measurable function $f$, $P_tf\in\Cn{1,2}(Q)$ for all $t>0$; see proof of \cite[Theorem 2.10.1]{Friedman69}.
\end{Remark}
\begin{Lemma}\label{lem:no_atom}
Let Assumption \ref{assp:regularity_Y} holds. Then, $\P(Y_T^{q,y}=0)=0$ for all $q\in\Qc$.
\end{Lemma}
\begin{proof} %By Assumption \ref{assp:pi_eta_lambda}-(iv), there exists a $\bar q>0$ such that for all $(t,y)$, $\beta\eta_t(q)-\gamma_t(y)\lambda_t(q)>0$ and $\tilde\pi_t(q)<0$ for $q\ge\bar q$. Therefore in problem \eqref{prob:control_problem} one can restrict to all strategies $q_\cdot\in\Qc$ bounded by $\bar q$. This, in particular, implies that
Consider the (not necessarily probability) measure $\tilde\P$ defined by
\[
\frac{d\tilde\P}{d\P}=\exp\left(-\int_0^T\zeta_tdW_t-\frac12\int_0^T\zeta^2_tdt\right)
\]
where $\zeta_t:=-\lambda_t(q_t)+\beta \eta_t(q_t)/\gamma_t(Y_t^q)$.  Then, one can write 
$\tilde\P(Y_T^{q,y}=0)=\E^{\P}[1_{\{\tilde Y_T^{y}=0\}}]$ where $\tilde Y$ under $\tilde \P$ satisfies
$d\tilde Y= \mu_t(\tilde Y_t)dt+\gamma_t(\tilde Y_t)dW_t$,
where $W$ is a Brownian motion under $\P$. Assumption \ref{assp:regularity_Y} implies that Aronson inequality holds for the density of $\tilde Y_T$; in particular, $\tilde Y_T$ has no atoms. Thus, $\tilde\P(Y_T^{q,y}=0)=0$ and since $\tilde\P<<\P$, $\P(Y_T^{q,y}=0)=0$. %{\Large\color{red}HERE I HAVE SERIOUS PROBLEMS USING GIRSANOV THEOREM AS NOVIKOV CONDITION FAILS TO HOLD FOR THE DENSITY.}
\end{proof}

\begin{Lemma}\label{lem:verification}
Let Assumption \ref{assp:regularity_Y} holds and triple $(\tilde\pi,\eta,\lambda)$ satisfy Assumption \ref{assp:pi_eta_lambda}. Suppose that
 $v\in\Cn{1,2,0}([0,T)\times\R\times\R)$ be such that $v_{e+}$ exists for all $(t,y,e)$. If $v$ is a supersolution of \eqref{prob:HJB_Appendix}, then $v\ge V$ for all $t<T$. In addition, if  there exists a measurable $q^*:=q^*(t,y,e)$ such that \eqref{eqn:Y^q} admits a strong solution and
\[
\begin{split}
 0
 &=
-\partial_t v -\frac{\gamma_t^2}{2}v_{yy}-\mu_tv_y-\tilde\pi_t(q^*)-\eta_t(q^*)v_{e+}-(\beta\eta_t(q^*)-\gamma_t\lambda_t(q^*))v_y \\
 v&(T,y,e)
 =
 -\alpha{1}_{\{y>0\}} e ,
 \end{split}
 \]
 Then, $V=v$.
\end{Lemma}
\begin{proof}
For the moment, let  $v\in\Cn{1,2,1}([0,T)\times\R\times\R)$. Then, for any $q\in\Qc$, It\^o's formula implies
\begin{equation}\label{eqn:ito_for_supersolution}
\begin{split}
v(\theta,E_\theta^{q},Y_\theta^q)&=v(t,y,e)+\int_t^\theta\bigl(\partial_t v +\frac{\gamma^2}{2}v_{yy}+\mu v_y+\eta_t(q)v_{e}\\
&+(\beta\eta_t(q)-\gamma\lambda_t(q))v_y\bigr)(s,E_s^{q},Y_s^q)ds+M_\theta-M_t
\end{split}
\end{equation}
where $M_\cdot$ is a continuous local martingale. Then, supersolution property of $v$ implies that
\[
v(\theta,E_\theta^{q},Y_\theta^q)\le v(t,y,e)-\int_t^\theta\tilde\pi_t(q)(s,E_s^{q},Y_s^q)ds+M_\theta-M_t.
\]
Let $\{\tau_n\}$ be a sequence for $\M_\cdot$  in the definition of  local martingale such that $\tau_n\to\infty$. By choosing $\theta=\tau_n\wedge T$, taking expectation $\E_{t,e,y}$, and sending $n\to\infty$, we obtain that 
\[
\begin{split}
v(t,y,e)&\ge \E_{t,e,y}\left[\int_t^\theta\tilde\pi_t(q)(s,E_s^{q},Y_s^q)ds+v(T,E_T^{q},Y_T^q)\right]\\
&=\E_{t,e,y}\left[\int_t^\theta\tilde\pi_t(q)(s,E_s^{q},Y_s^q)ds-\alpha{1}_{\{Y_T^{q}\ge0\}}E_t^{q}\right].
\end{split}
\]
The equality in the above holds from Lemma \ref{lem:no_atom}.
If $v\in\Cn{1,2,0}$, then by Krylov method of shaking coefficients \cite[proof of Theorem 2.2]{Krylov05}, one can find a supersolution $v_\eps(t,y,e)\in\Cn{1,2,1}$ such that $|v-v_\eps|=o(\eps)$. Thus, $v_\eps\ge V$ and the proof of the first part is complete after sending $\eps\to0$.

For the second part, one can see that all the above holds with equality, if one can show \eqref{eqn:ito_for_supersolution} holds with $v_{e+}$ in place of $v_e$. One can use a net of mollifiers $\{\rho_\eps\}_{\eps>0}$ with $\rho_\eps$ supported on $[0,\eps]$; i.e. $v_\eps:=v(t,\cdot,y)\ast\rho_\eps(e)$. Then one can write \eqref{eqn:ito_for_supersolution} for $v_\eps$. Since $\rho_\eps$ is supported on $[0,\eps]$, by sending $\eps\to0$, $\partial_ev_\eps\to v_{e+}$, which completes the proof.
\end{proof}

The first regularity result is covered by the following two lemmas.
\begin{Lemma}\label{lem:regularity_in_e}
$V$ is convex and continuous in $e$ uniformly on $(t,y)\in[0,T]\times\R$
\end{Lemma}
\begin{proof}
Convexity in $e$ follows from that $V$ is supremum of linear functions in $e$.
For  $q\in \Qc$, we can write 
\[
J_{q}(t,y,e)-J_{q}(t,y,e')=-\alpha\E\Bigl[1_{\{Y^q_T\ge0\}}\bigl(e-e')\bigr)\mid Y_t^q=y\Bigr],
\] 
where $E^q_T=\int_t^T\eta_s(q_s)ds$.
Thus,
$|J_{q}(t,y,e)-J_{q}(t,y,e')|\le \alpha|e-e'|$ and the inequality is uniform on $q\in\Qc$, which completes the proof. 
\end{proof}      
The following corollary follows from the properties of convex functions and the above Lemma.
\begin{Corollary}\label{cor:right_left_derivatives}
Right (left) partial derivatives of $V$, i.e. $\partial_{e+}V$  ($\partial_{e-}V$) exists, is non-decreasing and is right(left)-continuous  and bounded in $[-\alpha,0]$.
\end{Corollary}
\begin{Remark}\label{rem:equation_with_V_e+}
By Corollary \ref{cor:right_left_derivatives}, in Definition \ref{defn:viscosity} of viscosity supersolution solution, a continuously differentiable test function $\varphi$ which touches $V$ from below satisfies $\partial_{e-}V\le\partial_e\varphi\le\partial_{e+}V$. Therefore, the supersolution property implies that
\begin{align*}
-\frac{\partial \phi}{\partial t}-H^*(t,y,\phi_y,V_{e+},\phi_{yy})\ge 0
 \end{align*}
provided that $V$ is continuous for $t<T$.

In addition, if $H$ is continuous, then the above inequality must hold as equality. To see this, consider a point $(t_0,e_0,y_0)$ at which we have a strict inequality in the above. It follows from Lemma \ref{lem:regularity_in_e} and Corollary \ref{cor:right_left_derivatives} that there exists a point $e_1>e_0$ such that the above inequality holds at $(t_0,e_1,y_0)$ and $V_e$ exists at $(t_0,e_1,y_0)$. This violates the subsolution property of $V$ at $(t_0,e_1,y_0)$.
 
For the subsolution property, the set of the test functions $\varphi$ which touch $V$ from above is empty unless $V_{e+}=V_{e-}$, i.e. $V_{e}$ exists. 
\end{Remark}
Following the above remark, one can study the following terminal value problem which shares the same solution $V$ with  \eqref{prob:HJB_Appendix}.
 \begin{equation} \label{prob:HJB_with_V_e+}
 \begin{split}
 0
 &=
-\partial_t \bar V-H(t,y,\partial_y \bar V,\partial_{e+} V, \partial_{yy} \bar V) \\
 \bar V&(T,y,e)
 =
 -\alpha{1}_{\{y>0\}} e ,
 \end{split}
 \end{equation}
To establish the regularity property in $t$ and $y$, we present the following Lemma.
\begin{Lemma}\label{lem:regularity_in_t_y}
Let Assumption \ref{assp:regularity_Y} holds and triple $(\tilde\pi,\eta,\lambda)$ satisfy Assumption \ref{assp:pi_eta_lambda}. Then, for all $t<T$ and $e\in\R$, the partial derivatives $\partial_t V(t,y,e)$, $\partial_y V(t,y,e)$ and $\partial_{yy} V(t,y,e)$ exist and are continuous.
\end{Lemma}
\begin{proof}
First observe that by Assumption \eqref{assp:pi_eta_lambda}-(iv) and \cite[Theorem 1.7]{WangII92}, the viscosity solution $V$ to \eqref{prob:HJB_with_V_e+} is $\Cn{\frac12,1}$ in $(t,y)$ on $[0,T)\times\R$ for any fixed $e$. Now for fixed $e$, consider the following boundary value problem on parabolic domain $Q:=[t,t']\times[a,b]$
\[
\begin{split}
&0 = -\partial_t w-\mu w_y - \frac12\gamma^2w_{yy}-f(t,y),\; \text{\rm for}\; 0\le t<t',\;\text{\rm and }\; a\le y\le b\\ 
&w(t,a)=V(t,a,e),\; w(t,b)=V(t,b,e),\; \text{\rm for}\; 0\le t<t'\\
&w(t',y) =  V(t',y,e),\; \text{\rm for}\; a\le y\le b
\end{split}
\]
 where $-\infty<a<b<\infty$ and 
\[
f(t,y):=\sup_{0\le q\le\bar q} \left\{ \tilde\pi_t(q)+\eta_t(q)V_{e+}+(\beta\eta_t(q)-\gamma_t(\cdot)\lambda_t(q))V_y \right\}
\]
Since $f$ is locally bounded, one can apply Duhamel's principle on $Q$ to obtain
\[
w(t,\cdot)=P_{t'-t}w(T,\cdot)+a\int_t^{t'}P_{s-t}f(s,y)ds,\;\;\; \text{\rm for all }\;t<t\le T.
\]
By Assumption \ref{assp:regularity_Y}, the right hand side in the above and consquently $w$ is $\Cn{1,2}$ in $(t,y)$ for all $e$ and $t<T$. Notice that $\bar w:=V-w$ is a viscosity solution of 
$ 0=-\partial_t \bar w-\mu \bar w_y - \frac12\gamma^2\bar w_{yy}$ and $\bar w(T,y,e)=0$
 which has a uniques solution $\bar w\equiv 0$. Thus, $V$ is $\Cn{1,2}$ in $(t,y)$ for all $e$ and $t<T$.
%
%The opening argument in the proof of Lemma \ref{lem:no_atom} implies that in \eqref{prob:HJB_Appendix}, $H$ can be replaced by $H_{\bar q}$, i.e.
%\[ 
% 0
% =
%-\partial_t V-\mu V_y - \frac12\gamma^2V_{yy}-H_{\bar q}(t,y,V_y,V_e),\;\;\;\;\;
% V(T,y,e)
% =
% -\alpha{1}_{\{y>0\}} e.
%\]
\end{proof}

%\begin{Remark}\label{rem:regularity_linear}
%Let $\tilde\pi_t(q)=-\rho q^2+q$, $\eta_t(q)=\beta q$ and $\lambda_t(q)=\lambda_0+\lambda_1q$ for $\rho, \lambda_0,\lambda_1>0$ and assume that $\gamma_t(y)=\gamma$ is a positive constant. Then, triple $(\tilde\pi,\eta,\lambda)$ does not satisfy Assumption \ref{assp:pi_eta_lambda}-(iv) when $\beta<\gamma\lambda_1$. However, one can still show that the regularity result of Lemma \ref{lem:regularity_in_t_y} holds by applying the change of variable $W(t,y):=\exp(aV(t,y,e))$ for  some $a>0$. Here we consider variable $e$ as a fixed parameter. Then, $W$ satisfies
%$-\partial_t W -(\mu+\gamma \lambda_0) W_y - \frac{\gamma^2}{2}W_{yy}-gW=0$,
%where 
%\[
%g:=\frac{1}{4\rho}(1+V_e+(\beta-\gamma\lambda_1)V_y))_+^2-aV_y^2
%\]
%Then, appropriate choice of $a$ makes $g$ bounded, and the argument of Lemma \ref{lem:regularity_in_t_y} applies to show $W$ is $C^{1,2}$ in $(t,y)$.
%\end{Remark}
\begin{Lemma}\label{lem:optimal_control}
Let Assumption \ref{assp:regularity_Y} holds and triple $(\tilde\pi,\eta,\lambda)$ satisfy Assumption \ref{assp:pi_eta_lambda}. Then,  $q^*(t,y,e)$  given by 
\[
q^*\in\argmax_{q\ge0}\{\tilde\pi_t(q)-\eta_t(q)V_{e+}-(\beta\eta_t(q)-\gamma_t\lambda_t(q))V_y\}
\]
is an admissible Markovian  optimal control.
In addition, $q^*(t,y,e)$ is locally bounded, continuous in $(t,y)\in[0,T)\times\R$ and right-continuous in $e\in\R$.
\end{Lemma}
\begin{proof}
This is a direct consequence of Lemmas \ref{lem:regularity_in_e}, \ref{lem:regularity_in_t_y}, and \ref{lem:verification}.
\end{proof}
\begin{Lemma}\label{lem:V_e+=-S_t}
Let triple $(\tilde\pi,\eta,\lambda)$ satisfies Assumption \ref{assp:pi_eta_lambda} and Assumption \ref{assp:regularity_Y} holds true. Then, 
\[-V_{e+}(t,y,e)=S_t:=\alpha\P_{t,e,y}(Y_T^{q^*}\ge0),\]
where $q^*(t,y,e)$ is given by Lemma \ref{lem:optimal_control}.
\end{Lemma}
\begin{proof}
Suppose that $e'>e$. and let $q^*$ be the optimal control for problem \ref{prob:control_problem}  starting at $(t,y,e)$. 
Then, by direct calculations one can write
\[
V(t,y,e')-V(t,y,e)
\ge-(e'-e)\alpha\P_{t,e,y}(Y_T^{q^*}\ge0)
\]
Dividing both sides by $e-e'$ and sending $e'\to e$ yields to $V_{e+}(t,y,e)\ge -\alpha\P_{t,e,y}(Y_T^{q^*}\ge0)$.
One can obtain the other inequality by the fact that according to Lemma \ref{lem:optimal_control},  $q^*$ is right-continuous in $e$ and $Y_T^{q^*}$ has no atoms. If $q^*$ is the optimal control for problem \ref{prob:control_problem}  starting at $(t,y,e')$, then 
\[
V(t,y,e')-V(t,y,e)
\le-(e'-e)\alpha\P_{t,e',y}(Y_T^{q^*}\ge0).
\]
Sending $e'\to e$ yields to $V_{e+}(t,y,e)\le -\alpha\P_{t,e,y}(Y_T^{q^*}\ge0)$.
\end{proof}
 %%%%%%%%%%%%%%%%%%%%%%%%%%%%%%%%%%%%%%%%%%%%%%%%%%%%
%%%%%%%%%%%%%%%%%%%%%%%%%%%%%%%%%%%%%%%%%%%%%%%%%%%%
%%%%%%%%%%%%%%%%%%%%%%%%%%%%%%%%%%%%%%%%%%%%%%%%%%%%
\section{Numerical scheme}\label{appnedix_numerics}
In this section, we present details of numerical approximation of the nonlinear problem \eqref{eqn:pde_example} from Section \ref{Sec:numerics}. %It is worth mentioning that the techniques discussed here can be adjusted to a wider class of degenerate semilinear HJB equations.
The first step is to discretize in time and in $(y,e)$-space. Let $\Delta t:=\frac{T}{N}$ be the time step and $t^{(k)}=k\Delta t$, for $k=0,\cdots,N$. 
We set a computational bounded domain $[-L_e,L_e]\times[-L_y,L_y]$ for the $(y,e)$ space domain and discretize the computational domain by an appropriately fine grid $\{(e_i,y_j)\; : \; i=-N_e,...,N_e\; \text{ and } \; j=-N_y,...,N_y\}$ with $e_i=i\Delta e$, $y_j=j\Delta y$, $\Delta e=\frac{L_e}{N_e}$ and $\Delta y=\frac{L_y}{N_y}$. 
We set the discrete terminal data 
$V^N(e_i,y_j)=-\alpha{1}_{\{y_j\geq 0\}}e_i$.
To solve  \eqref{eqn:pde_example} numerically, we need to set (1)
 appropriate artificial boundary conditions (a.k.a. ABC) for the computational domain, 
(2) treatment of  discontinuity of the terminal condition, and (3) stable approximation of the semi-linear terms in \eqref{eqn:pde_example}.

To properly set the ABC for computational domain, we return to the optimization problem \eqref{eqn:large_lambda_impact}. If $L_y$ is sufficiently large so that $Y^q$ defined by \eqref{eqn:Y_under_Q^q} satisfies $\Q^q(Y^q_T\ge0\mid Y^q_t=L_y)\approx1$ uniformly on $q\in\Qc$, then we can approximate the value function with the following simple deterministic control problem.
\[
\begin{split}
V(t,e,L_y)&\approx\sup_{q_\cdot\in\Qc} \Bigl\{\int_t^T(-\varrho q_s^2+q_s)ds
                                -\alpha\Bigl(e+\int_t^Tq_sds\Bigr)
                          \Bigr\}\\
                          &=-\alpha e+\frac{(1-\alpha)^2(T-t)}{4\rho}
\end{split}
\]
 In a similar but simpler fashion, at $y=-L_y$ we have $\Q^q(Y^q_T\ge0\mid Y^q_t=-L_y)\approx0$, and thus the approximate ABC becomes $V(t,e,-L_y)\approx\frac{1}{4\varrho}(T-t)$. We postpone the derivation of boundary condition on $e=L_e$ or/and $e=L_e$ for after we present the algorithm.

In order to handle the discontinuity of terminal condition in the algorithm, we adopt a splitting (Trotter-Kato type) method. At each time step, we handle the calculations in to half-steps. In the first half-step, we solve the heat equation $-\partial_tv -\frac{\gamma^2}{2}\partial_{yy}v=0$ with the same boundary conditions as in the previous step. This  regularizes the discontinuous terminal condition. In the second half-step, we solve $-\partial_tv-\mu\partial_yv-(1+\partial_ev+(1-\gamma)\partial_y v)^2_+/4\rho=0$ with ABC boundary condition for the current time step.

To treat the semi-linear term $(1+\partial_eV+(1-\gamma)\partial_y V)^2_+/4\rho$ in the second half-step above, we write it as the multiplication of $(1+\partial_eV+(1-\gamma)\partial_y V)_+/2\rho$ and $(1+\partial_eV+(1-\gamma)\partial_y V)/2$. Notice that the first term is equal to the optimal control $q^{(3)}$. If we calculated the first term by using the first half-step (solution heat equation), then the second half-step is to solve a linear equation $-\partial_tv-\mu\partial_yv-q^{(3)}(1+\partial_ev+(1-\gamma)\partial_y v)/2=0$.
The above discussion is summarized in the following algorithm.

\begin{algorithm}
\caption*{The splitting scheme for problem \eqref{eqn:pde_example}}
\begin{algorithmic}[1]
%\Procedure{} {}
\\ $\hat V^{N}(T,e_i,y_j)=-\alpha 1_{\{y_j\ge0\}}e_i$.
\For{each $n=N-1,...,0$ }
\\$\hat V^{n+\frac12}(y,e):=V(t_n,y,e)$ where $V$ is the solution of $\partial_tV+\frac{1}{2}\gamma^2V_{yy}=0$ on $[t^{(n)},t^{(n+1)}]$ with boundary condition $V(t,e,-L_y)=\frac{(T-t_{n+1})}{2\varrho}$ and  $V(t,e,-L_y)=U(t_{n+1},e)$  and terminal condition $V(t_{n+1},y,e)=\hat V^{n+1}(y,e)$.
\\ $\varphi^{n}(y,e):=\frac{1}{2\varrho}(1+\hat V^{n+\frac12}_e+(1-\gamma)V^{n+\frac12}_y)_+$.
\\$\hat V^{n}(y,e):=V(t_n,y,e)$ where $V$  is the solution of 
\begin{equation}\label{eqn:1st_order_pde}
\partial_tV+\mu V_y+\frac{\varphi^n}{2}\left(1+V_e+(1-\gamma)V_y \right)=0
\end{equation}
 on $[t^{(n)},t^{(n+1)}]$ with boundary condition $V(t,e,-L_y)=V(t,-L_y,e)=\frac{(T-t_{n})}{2\varrho}$,  $V(t,e,-L_y)=U(t_n,e)$,  and  terminal condition $V(t_{n+1},y,e)=\hat V^{n+\frac12}(y,e)$.
\EndFor
%\EndProcedure
\end{algorithmic}
\end{algorithm}
To avoid the hassle of setting ABC on both $e=L_e$ and $e=-L_e$, we can approximate $V_e(t_n,e_i,y_j)$ from one side by $\frac{V(t_n,e_{i+1},y_j)-V(t_n,e_i,y_j)}{\Delta e}$. 
To set the boundary condition at $-L_e$, notice that first order linear PDE \eqref{eqn:1st_order_pde} can easily be solved by the method of characteristics. However, we can only use method of characteristics as long as we stay in the computational domain. More specifically, method of characteristics can give us the solution at $e=-L_e$; i.e.
\begin{equation}\label{eqn:char}
\begin{split}
\hat V^{n}(t_n,-L_y,e)&=\hat V^{n+\frac12}(e+\psi^n(-L_y,e),y+(1-\gamma)\psi^n(-L_y,e)+\mu\Delta t)\\
&\hspace*{4cm}-\varrho (\psi^n(-L_y,e))^2+\psi^n(-L_y,e).
\end{split}
\end{equation}
Therefore, we use eqref{eqn:char} to set ABC at $e=-L_e$ for \eqref{eqn:1st_order_pde} at  step 5 in the algorithm is solved, i.e. 
\[
V(t,0,y)=\hat V^{n+\frac12}(e+\psi^n(0,y),y+(1-\gamma)\psi^n(0,y)+\mu\Delta t)-\varrho (\psi^n(0,y))^2+\psi^n(0,y).
\]
\begin{Remark}\label{rem:char=dpp}
Estimation \eqref{eqn:char}, which is based on the method of characteristics, can be equivalently derived from approximate dynamic programing principle for the following deterministic optimal control problem which corresponds to \eqref{eqn:1st_order_pde}.
\begin{equation}\label{eqn:deterministic_optimal_control}
\sup_{q\in\Qc}\int_0^T(-\varrho q_t^2+q_t)dt-\alpha E^q_T1_{\{Y_T^q\ge0\}},
\end{equation}
where $dE_t^q=q_tdt$ and $dY_t^q=(\mu+(1-\gamma)q_t)dt$.
 The dynamic programming principle of problem \eqref{eqn:deterministic_optimal_control} over the interval $[t^{(n)},t^{(n+1)}]$ is 
\[
\begin{split}
\ V(t_n,0,y)&=\sup_{q\in\Qc} \int_{t_n}^{t_{n+1}}(-\varrho q_t^2+q_t)dt\\
&+V\Biggl(t_{n+1},e+\int_{t_n}^{t_{n+1}}q_tdt,y+(1-\gamma)\int_{t_n}^{t_{n+1}}q_tdt+\mu\Delta t\Biggr)
\end{split}
\]
Observe that $\psi^n(0,y)$ is an approximation of the optimal control on interval $[t^{(n)},t^{(n+1)}]$. Thus replacing $q_t$ by $\psi^n(0,y)$ yields \eqref{eqn:char}.
\end{Remark}

\bibliographystyle{plain}
\bibliography{Carbon-2015-09-04.bbl}

\end{document}